# Causal Inference Through Potential Outcomes and Principal Stratification: Application to Studies with "Censoring" Due to Death[1]

**Donald B. Rubin**

*Abstract.* Causal inference is best understood using potential outcomes. This use is particularly important in more complex settings, that is, observational studies or randomized experiments with complications such as noncompliance. The topic of this lecture, the issue of estimating the causal effect of a treatment on a primary outcome that is "censored" by death, is another such complication. For example, suppose that we wish to estimate the effect of a new drug on Quality of Life (QOL) in a randomized experiment, where some of the patients die before the time designated for their QOL to be assessed. Another example with the same structure occurs with the evaluation of an educational program designed to increase final test scores, which are not defined for those who drop out of school before taking the test. A further application is to studies of the effect of job-training programs on wages, where wages are only defined for those who are employed. The analysis of examples like these is greatly clarified using potential outcomes to define causal effects, followed by principal stratification on the intermediated outcomes (e.g., survival).

## PROLOGUE

This article is the written version of two presentations that I had the privilege of giving in 2005, the first at Carnegie Mellon University on September 16, the Morris DeGroot Memorial Lecture, and the second at the Washington, D.C., Chapter of the American Statistical Association's Morris Hansen Memorial Lecture, November 2. Both were truly enjoyable and stimulating occasions for me, not only the presentations themselves, but the warm events following them. I am extremely grateful to the two selection committees for inviting me, to the relatives of both Morrises, and to the very good friends who shared time with me in Pittsburgh and Washington.

The basic material in the talk had been presented a couple of previous times before the first "Morris" talk in Pittsburgh, so it was fairly polished, I thought, and thus worthy of memorial lectures honoring these two wonderful statisticians and good

*Donald B. Rubin is John L. Loeb Professor of Statistics, Department of Statistics, Harvard University, Cambridge, Massachusetts 02138, USA e-mail: rubin@stat.harvard.edu.*

friends. Also, the content was rather broadly accessible and relevant to both men's interests.

I do not know whether Morrie (DeGroot) or Morris (Hansen) knew each other well or not; they tended to travel in different statistical circles, Morrie more in the Bayesian decision theory, academic statistics world, and Morris more in the survey design, government statistics world, but they were both very influential and widely admired.

A few words about Morrie first, partly because the Morrie talk was first, but also because I met Morrie first, although I was certainly familiar with both names as a graduate student at Harvard in the late 1960s. In 1976 Morrie became Editor of *Journal of the American Statistical Association*, Theory and Methods, and he contacted me to stay on as an Associate Editor—of course, I was thrilled and agreed. But I wondered why because I didn't really know him at all. He explained that when he had been a *JASA* Associate Editor earlier, he had sent me various papers to review, and he liked the reports that I wrote. In particular, one sequence of papers that he sent me to review was on a topic that I had felt was completely old-fashioned and unimportant, and my reviews repeatedly said so. I had asked Morrie at that time why he kept sending them to me since I didn't like them. And he replied that besides me, there appeared to be only two kinds of possible reviewers: the vast majority who refused to read the submissions because they were negative about the area and didn't want to waste their own time; and a very few who loved the boring material because it was what they did and therefore would uncritically recommend publication. Morrie wanted ammunition to recommend rejection, which would be difficult with only positive reviews, and I was providing that ammunition! As Editor, he promised to use me more efficiently except when a paper on this particular topic arrived.

There are many other Morrie stories available, some from the very early Bayesian meetings in Valencia, a quarter of a century ago. One has Morrie explaining at his after-dinner talk how he could manage to stay up partying every night until the wee hours of the morning, consuming alcohol and smoking cigars, and still arise for an early breakfast, more jovial and energetic than those half his age—he explained that it was simple: practice, practice, practice. Good advice I've tried to follow.

My experience with Morris was more limited, and involved discussions and meetings often stimulated by survey nonresponse issues, or by his advisory roles at the Census, or on government committees. He was always warm but principled, with a keen desire to see statistics used to address important real-world problems. In many ways, he reminded me of my wonderful Ph.D. advisor, Bill Cochran, with respect to having similar attitudes toward the field of statistics and the problems it should be addressing.

I think that the topic of today's talk, and I hope, the presentation itself, would be of interest to both of these pillars of twentieth century statistics.

## 1. INTRODUCTION TO DATA THAT ARE "CENSORED" OR "TRUNCATED" DUE TO DEATH

There are several themes in this presentation that are quite general. First, the proper analysis of complicated randomized experiments can often take on many of the features of the proper analysis of non-randomized (observational) data; in both, covariates play an important role, which is often unappreciated. Second, it is critical to give adequate conceptual thought to any nonstandard statistical problem before attacking it with mathematical analysis or available computer programs. As Picasso said: "computers are worthless; they only give answers." (Thanks to Stuart Baker for first pointing out this great quote.) And third, intermediate outcome variables, which arise frequently in practice but often remain unrecognized, are not easy to handle well; in fact, the giant of statistics, Sir Ronald Fisher, gave flawed advice about them throughout his career (Rubin, 2005). To be fair, however, nearly all researchers I have read have also failed to provide good advice on this tricky topic of intermediate outcome variables, and Fisher appeared never to have focused any real attention on it.

One generic example of a complicated randomized experiment with an intermediate outcome variable is the specific topic of this presentation, and can be labeled as involving "censoring" or "truncation" of data due to death. For instance, the patient in the experiment dies after treatment assignment, but before the primary outcome variable, say Quality of Life (QOL) two years after assignment, can be measured. An artificial example of this will be used throughout this presentation.

Examples of such censoring also exist in other fields. For instance, suppose that we were interested in the effect of a special educational intervention in high school on final test scores, and some of the students in a randomized experiment evaluating this intervention do not finish high school. Or in some economics situations, interest focuses on the causal effect of a job-training program on wages (not income), which are only well defined for those people who are employed; thus, people who are unemployed when wages are measured have their wage outcome data "censored" or "truncated." Or, suppose in a study of the effects of hormone replacement therapy (HRT) on five-year cancer-free survival, some women die before five years, but are cancer-free when they die, say, of heart disease at three years. As this short list of examples makes clear, this type of complication can and does arise in many circumstances.

My first contact with this specific issue was in the context of a consulting project in the early 1990s for AMGEN for a product for the treatment of ALS (amyotrophic lateral sclerosis) or "Lou Gehrig's Disease"—see a brief discussion in Rubin (2000), and prior to that in Rubin (1998). ALS is a progressive neuromuscular disease that eventually destroys motor neurons, and death follows, typically from lungs that are unable to operate. No good treatments were (or are) available. In the AMGEN example, the active treatment, say product T, was to be compared to the control treatment, C, where the primary outcome was QOL two years post-randomization, as measured by "forced vital capacity" (FVC), essentially, how big a balloon you can blow up when you are alive. When FVC is large, you can typically get on fairly well, whereas when this is small, you are in very bad shape. In fact, many people do not reach the end-point of two-year post-randomization survival, and so two-year QOL is "truncated" or "censored" by death. I was brought into this project because, as sometimes is the case, the unavailable QOL data were trying to be fit into a "missing data" framework.

Before continuing with this example, it is helpful to state that the general attack on this problem being presented here uses the framework of "principal stratification" (Frangakis and Rubin, 2002). The specific technical work on this topic was initiated in a Ph.D. thesis at Harvard University (Zhang, 2002), and follow-up work appears in Zhang and Rubin (2003) and Zhang, Rubin and Mealli (2005, 2006). These references provide discussion of other techniques that have been proposed to attack this problem, and why those techniques are generally deficient relative to the principal stratification approach presented here. We only briefly review these other deficient approaches later, after setting up a correct framework.

The key idea of principal stratification is to stratify on the intermediate outcome, here the indicators for two-year survival, but not on the *observed* two-year survival, which is an outcome generally affected by the treatment received. Rather we should stratify on the *bivariate outcome*: survival if assigned active treatment, survival if assigned control treatment. This bivariate outcome is not affected by the treatment received, even though which of the two outcomes is actually observed is affected by the treatment received. Thus, in our running example there are four principal strata representing four types of people: those who will live no matter how treated ($LL$), those who will die no matter how treated ($DD$), those who will live if treated but die if not treated ($LD$), and those who will die if treated but live otherwise ($DL$). A specific artificial case is displayed in Table 1 and will be used for most of this article. It is chosen to be relatively extreme to make points more dramatically; it does not realistically represent any data from the AMGEN trial, which originally motivated this approach.

## 2. THE RUNNING EXAMPLE

Table 1 presents the hypothetical truth, and displays what would happen to the group of people in each principal stratum under both the active treatment and the control treatment. Of course, for any person we can only observe the "potential outcomes" under one or the other treatment, not both—the fundamental problem facing (Rubin, 1978, § 2.4; Holland, 1986, §3). Holland called the general perspective to causal inference presented here the "Rubin Causal Model" (RCM) for a series of papers written in the 1970s expounding and expanding this perspective (Rubin, 1974, 1975, 1976, 1977, 1978, 1979, 1980); Table 1 assumes "SUTVA" (Rubin, 1980, 1990), the stable-unit-treatment-value assumption, or stability; this assumption is very commonly made.

TABLE 1
*Principal strata among the patients*

| % population | Principal stratum | Treatment | | Control | | Treatment effect on QOL |
|---|---|---|---|---|---|---|
| | | $S_i(T)$ | $\bar{Y}_i(T)$ | $S_i(C)$ | $\bar{Y}_i(C)$ | |
| 20 | $LL$ | $L$ | 900 | $L$ | 700 | 200 |
| 40 | $LD$ | $L$ | 600 | $D$ | * | * |
| 20 | $DL$ | $D$ | * | $L$ | 800 | * |
| 20 | $DD$ | $D$ | * | $D$ | * | * |

The first row of Table 1 shows that 20% of the population will live under either treatment, as indicated by the survival potential outcomes $S(T) = L$ and $S(C) = L$; $S(T)$ is the potential outcome for survival when assigned treatment and $S(C)$ is the potential outcome for survival when assigned control. For these $LL$ people, the average $Y$ (i.e., QOL) if all were treated would be 900, which is good, but would be 700 if not treated, which is fair. Therefore, the average causal effect of the treatment for the $LL$ stratum is $900 - 700 = 200$, as indicated in the last column. This will be called the SACE—the survivor average causal effect. Critically, a causal effect must be a comparison of treatment potential outcomes, $Y(T)$, and control potential outcomes, $Y(C)$, on a *COMMON* subset of units, here the set of $LL$ units.

The $LD$ units, displayed in the second row of Table 1, are those who would die under control but live under treatment [i.e., $S(T) = L$ and $S(C) = D$], and they comprise 40% of the population. If these units were all treated, their average QOL would be 600, which is poor, but if they were not treated, they would die, and their QOL would be undefined (or defined on the sample space of the positive real numbers extended to include an asterisk). To assign a particular value to QOL when dead is to assume we know how to trade off a particular QOL and being dead (and out of misery). Not only do we not know how to do this, but the trade-off could vary by individual, so we prefer simply to represent the actual truth at this point, and not bring in such extraneous value judgments.

The third row of Table 1 is for those who would die under treatment but live under control, those in the $DL$ group with $S(T) = D$ and $S(C) = L$. These subjects comprise 20% of the population, and their average QOL under control is a quite decent 800. And the final 20% represented in the fourth row are in the $DD$ group, who would die no matter which treatment they received.

A well-defined real value for the average causal effect of the active treatment versus the control treatment on QOL exists only for the $LL$ group. For the $LD$ and $DL$ groups, the average causal effect on QOL involves the aforementioned trade-offs with death, and for the $DD$ group there is no QOL to compare, so the causal effect on QOL for them must be zero. The most that we can ever hope to learn in any study of this population of values under these two treatments is recover this table of values.

Before considering how to do this, however, let us examine this table a bit more. First, the active treatment is better for survival than the control treatment because 60% (20% $LL$ + 40% $LD$) would survive when treated, whereas only 40% (20% $LL$ + 20% $DL$) would survive if not treated (control).

Thus the active treatment is better for overall survival, and the active treatment is better for QOL for the subset of people where it is well defined, the $LL$ group, by $+200$. Therefore, with no more information about possible subgroup differences, such as differences between males versus females, the treatment is preferable for the population.

Notice also in this example that even if all four groups had been the same size, each representing a quarter of the population, the treatment still would have been preferred to control. The reason is that, although there would have been no treatment versus control difference on overall survival, treatment would have a positive causal effect on QOL for the only subgroup where it is well defined. If, in this case, an * were imputed with 0, the conclusion would have been that there is no benefit to the active treatment for either survival or QOL because the last column would have averaged to zero $(200 + 600 - 800 + 0)/4 = 0$, a conclusion that conflates facts with value judgments. This conclusion would be especially deceptive if conclusions from this population were to be generalized to future healthier populations dominated by people like those in the $LL$ group; this often can occur in real-world clinical trials, where experimental drugs are first tried with sicker patients, and approval is based on results with these patients, but if approved, the drugs are used in broader and healthier populations.

Continuing with the examination of Table 1, under treatment the healthiest group is the $LL$ group, followed by the $LD$ group; the $DL$ and $DD$ groups both die when treated. However, under control the *DL* group is healthier than the *LL* group, and both of these groups are healthier than the *LD* and *DD* groups, whose members die under control. Can this be realistic? The answer is "yes" for at least two reasons: First, some drugs do have negative side effects for some subgroups of people, and so here that would be the *DL* subgroup, who would survive if untreated. Another possible reason is that the active treatment may make some people feel so much better, even though it does not affect their disease progression, that they "overdo" it—play tennis, go to parties, have normal sex lives, and so on. There are some drugs that can have effects like this; Epogen, another product made by AMGEN, substantially increases red blood cell production and is of substantial apparent benefit to dialysis and chemotherapy patients, who can have much more energy with the extra oxygen-carrying capability created by Epogen. For example, Epogen has become an issue in recent years in some professional sports (e.g., bicycling with Lance Armstrong recently, and Jerome Chiotti before him). These situations reinforce the related points made earlier about the trade-offs between a potentially higher quality of life versus an earlier death. For example, a weak 90-year-old may consider a QOL of 600 preferable to death, whereas an Olympic athlete who is used to running ten miles a day may prefer death to a completely sedentary and deteriorating QOL.

A related point is that Table 1 is only a summary of the individuals' potential outcomes in this hypothetical population because it only gives the mean values of the survival and QOL potential outcomes within each principal stratum. The more complete version of Table 1 would also provide the marginal distributions of all four potential outcomes, in addition to their means, and moreover, would provide the joint four-dimensional distribution of the potential outcomes within each principal stratum. Having such information would allow individuals to make the trade-off between death and QOL, but it is far more difficult to estimate such a table of joint distributions than simply the means, because treatment and control potential outcomes are never jointly observed. We return briefly to this topic after understanding the simpler problem of estimating the means given in Table 1 from observable data.

For now let us accept Table 1 with just $Y$ means in all four principal strata as truth, and consider next how we learn about this table from observable data.

## 3. WHAT WOULD BE OBSERVED IN A RANDOMIZED EXPERIMENT?

Suppose that we conducted a huge completely randomized experiment on a huge random sample from this population: half get randomized to active treatment and half get randomized to the control. Even though not blocked on the unknown principal strata, in expectation, half of each principal stratum will be exposed to each treatment. This is reflected in Table 2 where each row in Table 1 is split in half, with the top one in each half getting the active treatment, indicated by $Z = 1$, and the other getting the control treatment, indicated by $Z = 0$. Of course, we do not get to observe all the values in Table 2, and in fact, do not know the principal strata to which individual people actually belong.

Suppose now that we permute the rows in Table 2 so that rows that have the same observed treatment and the same observed survival are adjacent. We trivially obtain Table 3, but of course we still do not know the splits between the pairs of adjacent rows. That is, for the first pair of rows, we observe that all these are people who got treated and lived, and that these comprise 30% of the people in the experiment; consequently, the observed survival rate in the random half assigned treatment is 60% (= 2 × 30%). The average QOL for this group will be a $1/3 + 2/3$ mixture of $LL$ and $LD$, that is, of the averages 900 and 600, and so the observed average for those who got assigned treatment and lived will equal 700. We do not observe any control potential outcomes for these people because they are all treated. For the second pair of rows in Table 3, we have that they are observed to be treated and die, and comprise 20% of the population, or 40% of the treated group dies. Again, no control potential outcomes are observed for these people.

TABLE 2
*Principal strata among the patients, each split by treatment assignment*

| % population | Principal stratum | Assignment $Z_i$ | Treatment | | Control | | Treatment effect on QOL |
|---|---|---|---|---|---|---|---|
| | | | $S_i(T)$ | $\bar{Y}_i(T)$ | $S_i(C)$ | $\bar{Y}_i(C)$ | |
| 10 | LL | T | L | 900 | L | 700 | 200 |
| 10 | LL | C | L | 900 | L | 700 | 200 |
| 20 | LD | T | L | 600 | D | * | * |
| 20 | LD | C | L | 600 | D | * | * |
| 10 | DL | T | D | * | L | 800 | * |
| 10 | DL | C | D | * | L | 800 | * |
| 10 | DD | T | D | * | D | * | * |
| 10 | DD | C | D | * | D | * | * |

For the third pair of rows, we observe that they are assigned to control and live, and comprise 20% of the population, implying a survival rate in the population under control of 40%. Also they have an observed average QOL of 750, which arises from the $1/2 + 1/2$ mixture of $LL$ and $DL$ with means 700 and 800, respectively. In these two rows the treated potential outcomes are not observed because the people were assigned control. Finally, the last pair of rows were assigned control, and they were observed to die, and they comprise 30% of the population, implying a death rate in the control group of 60%. No treated potential outcomes are observed for these people.

The discussion in the previous two paragraphs is summarized in Table 4, which displays only what is actually observed in the study. Several features are noteworthy. First, suppose that we decide to assess the causal effect of the active control treatment on survival, the "intermediate" outcome. We get the correct answer: 60% survive when treated versus 40% when untreated, just as in Table 1. Next suppose that we decide to assess the causal effect of active treatment versus control treatment on QOL using only the subjects for whom we have observed QOL: we would compare the observed average of 700 for the treated group versus the observed average of 750 for the control group and conclude that, although the active treatment is good for survival, it is bad for QOL if you do survive—but this is simply wrong! The causal effect of the active treatment versus the control is positive (+200) for the $LL$ group, which is the only group for which QOL is well defined.

What went wrong with this last analysis comparing mean QOL for survivors? The answer is that the comparison does not estimate a causal effect. Rather than comparing treated and control potential outcomes [i.e., $Y(1)$ and $Y(0)$] on a common subset of units (like the $LL$ group), it compares the average observed treatment potential outcome $Y(1)$ of 700, which comes from a $1/3 + 2/3$ mixture of $LL$ and $LD$, with the average observed control potential outcome $Y(0)$ of 750, which comes from a $1/2 + 1/2$ mixture of $LL$ and $DL$. These are different groups of people, having only some $LL$ people in common, but even these are in different fractions.

This method of attack, comparing QOL when it is observed and dropping people who died, although popular in some settings, is simply wrong in general. But then what should we do? If we knew the labels of the principal strata for all the people, we could simply analyze the data within each stratum, in particular compare $Y(1)$ and $Y(0)$ in the $LL$ stratum, but we do not have this information. As Table 5 displays, instead, each of our observed groups defined by observed treatment assignment $Z$ and observed survival $S^{\text{obs}}$, comprises a mixture of people from two unobserved, or latent, principal strata.

TABLE 4
*Observed data for the example of Tables 1–3*

| % population | $Z_i$ | Treatment | | Control | |
|---|---|---|---|---|---|
| | | $S_i(T)$ | $\bar{Y}_i(T)$ | $S_i(C)$ | $\bar{Y}_i(C)$ |
| 30 | T | L | 700 | ? | ? |
| 20 | T | D | * | ? | ? |
| 20 | C | ? | ? | L | 750 |
| 30 | C | ? | ? | D | * |

TABLE 5
*Group classification based on observed treatment assignment and observed survival indicator $OBS(Z, S^{\text{obs}})$, and associated data pattern and possible latent principal strata*

| Observed group $OBS(Z, S^{\text{obs}})$ | $Z$ | $S^{\text{obs}}$ | $Y^{\text{obs}}$ | Possible latent principal strata |
|---|---|---|---|---|
| $OBS(T, L)$ | T | L | $\in R$ | $LL$, $LD$ |
| $OBS(T, D)$ | T | D | * | $DL$, $DD$ |
| $OBS(C, L)$ | C | L | $\in R$ | $LL$, $DL$ |
| $OBS(C, D)$ | C | D | * | $LD$, $DD$ |

TABLE 3
*Permuted table for principal strata among the patients, each split by treatment assignment*

| % population | Principal stratum | $Z_i$ | Treatment $S_i(T)$ | Treatment $\bar{Y}_i(T)$ | Control $S_i(C)$ | Control $\bar{Y}_i(C)$ | Treatment effect on QOL |
|---|---|---|---|---|---|---|---|
| 10 | LL | T | L | 900 | L | 700 | 200 |
| 20 | LD | T | L | 600 | D | * | * |
| 10 | DL | T | D | * | L | 800 | * |
| 10 | DD | T | D | * | D | * | * |
| 10 | LL | C | L | 900 | L | 700 | 200 |
| 10 | DL | C | D | * | L | 800 | * |
| 20 | LD | C | L | 600 | D | * | * |
| 10 | DD | C | D | * | D | * | * |

## 4. POSSIBLE APPROACHES

One possible approach is to treat the problem as one of missing data, and try to impute, or multiply impute, the "missing data" that are "censored" by death. But we really already rejected this idea in the discussion of Section 2: the $Y$ outcomes are not missing; they are undefined, or defined to be $*$. Maybe adding some simplifying assumptions would help?

There are some assumptions that are relatively standard in similar settings, in particular, where the intermediate outcome variable indicates compliance with assigned drug (Angrist, Imbens and Rubin, 1996). In Tables 1–4, the four principal strata could be called compliers, never-takers, always-takers and defiers, corresponding to $LD$, $DD$, $LL$ and $DL$, respectively. The analogy here is that people who would live under treatment but would die under control are "complying" with the encouragement of the active treatment to help them.

The first standard assumption (after SUTVA) that is made in the noncompliance setting is called "monotonicity" or the "no-defier" assumption, which rules out the $DL$ group. The no-defier assumption can be reasonable in our QOL setting, but is wrong in the context of our numerical example because there exist both $LD$ and $DL$ groups.

The next assumption that is often made in the noncompliance setting is called "exclusion," which asserts that if treatment assignment cannot change the intermediate outcome, $D$, it cannot change the final outcome, $Y$; here, that would mean that there is no treatment effect on $Y =$ QOL for either the $LL$ group or the $DD$ group. But the causal effect on QOL for the $LL$ group is precisely what we want to estimate, and so we cannot assume it to be zero! So the exclusion assumption does us absolutely no good. When both monotonicity and exclusion hold, however, the classical instrumental variables estimate (IVE) can be used to estimate the "complier average causal effect"; see Angrist, Imbens and Rubin (1996) for extended discussion. In simple settings, the IVE is the simple treatment minus control estimate for the mean of $Y$ divided by the simple treatment minus control estimate for the mean of $D$, so here would require some imputation of the "missing" $Y$ values for those who are observed to die. If we imputed zero for the QOL for those who die, we would have that IVE $= (420 - 300)/(0.6 - 0.4) = 600$, unrelated to anything real, which is not surprising because the underlying assumptions justifying the IVE are both wrong in our example.

Another possible assumption, considered in Zhang and Rubin (2003), is "stochastic dominance," which implies that, on average, the $LL$ group is healthier under control than the $DL$ group, and the $LL$ group

TABLE 6
*Bounds for treatment effect on QOL in LL for numerical example*

| Monotonicity assumption | Stochastic dominance | [Lower Bound, Upper Bound] |
|---|---|---|
| No | No | [−200, 200] |
| Yes | No | [−150, 0] |
| No | Yes | [−100, 150] |
| Yes | Yes | [−50, 0] |

The first column shows whether the monotonicity assumption (A1) is made, and the second column shows whether the stochastic dominance assumption (A2) is made. The last column shows the bounds for the numerical example of Tables 1–4.

is healthier under treatment than the *LD* group. Again, this condition is violated in our numerical example, as was noted in Section 2 earlier. Large-sample bounds on the Survivor Average Causal Effect (i.e., the causal effect in the *LL* group) are derived in Zhang and Rubin (2003) under monotonicity and stochastic dominance, but are not very useful in our example, as displayed in Table 6, because none of the assumptions holds. In other examples, they could be quite useful.

## 5. THE ROLE FOR COVARIATES

A more successful general approach is to collect and use covariates that are predictive of both the intermediate potential outcomes (e.g., here survival) and the final potential outcomes (e.g., here QOL). This was done in the actual AMGEN application because there were several measurements of baseline FVC (= baseline QOL). Thus at baseline measurements of each patient's current FVC and the rate of deteriorating FVC were available, and these were highly predictive of both survival and of two-year-later QOL if surviving.

To amplify this point, consider Table 7, which is identical to Table 1 except with an added left-most column, labeled $X$, for a covariate, baseline QOL in the hypothetical example. The hypothetical means of $X$ are displayed, and we will assume that the hypothetical variances of $X$ within each principal stratum are small relative to the differences between the means. We again pretend that we conduct a huge randomized experiment, with 50% treated and 50% control, to obtain Table 8, which parallels Table 2 with treated and control pairs of rows. Trying to permute the rows in this table to bring groups adjacent that are observed to be the same (with respect to observed $X$, observed treatment assigned, and observed survival) leaves the table unchanged because, for example, although the first and third rows are both treated and survive, they are distinguishable by their differing baseline QOL distributions. The observed data then are as in Table 9, from which we reach the following conclusions.

People with baseline FVC around 800, which is pretty good, comprise 20% of the population, and they all survive no matter how treated; the causal effect on QOL of active versus control for them is $900 - 700 = +200$—this conclusion agrees with the truth in Tables 1 and 7. Next, consider the 40% of the population with baseline QOL around 500, which is quite poor. They will survive if treated, with an average QOL of 600, better than at baseline, but still poor; without the active treatment, however, they will die; again, some may prefer death to poor QOL. For the 20% with baseline FVC around 300, which is very poor, neither active treatment nor control can prevent death. And finally, for the 20% with the best baseline FVC, around 900, which is quite good, we see that if not actively treated, their QOL will decline to 800, still not bad, whereas if treated they will die. This is unexpected and requires follow-up interviews with their individual doctors (and/or friends and spouses) to determine the reasons for their deaths, possibly negative side effects of the drug, or overactivity due to the drugs' dramatic effects on perceived health. But the point is that a "super" covariate has allowed the recovery of the truth in Table 1.

Notice that collecting more measurements of outcomes does not help in the same way as collecting covariates, because outcome measurements have different potential outcomes depending on the treatment exposures, and so using outcomes to improve inference will require some serious modeling efforts involving new assumptions. Some work on this topic appears in Zhang (2002), but is an important area for statistical research because it is common, and often easy, to collect such repeated measurements post-randomization.

Another source of information that can be utilized is the distributional shape of the outcome in the different groups and treatments. For example, if we knew that QOL measurements were approximately normally distributed across subjects within each principal stratum and treatment condition, this could be very helpful because it would allow standard mixture modeling tools to be used to help disentangle the normal components (e.g., Dempster, Laird and Rubin, 1977; Titterington, Smith and Makov, 1985). This approach is used in the example in Section 6.

## 6. THE ROLE FOR DISTRIBUTIONAL ASSUMPTIONS, SUCH AS NORMALITY

Here we extend the example in Table 1 to include the distribution of QOL within each of the four principal strata to illustrate how such information can be used to help recover the information in Table 1. This extension will also lead to a brief discussion of the more difficult issue of the role of the joint

Table 7 (= Table 1 with Key Covariate)
*Principal strata among the patients*

| $X$ | % population | Principal stratum | Treatment $S_i(T)$ | Treatment $\bar{Y}_i(T)$ | Control $S_i(C)$ | Control $\bar{Y}_i(C)$ | Treatment effect on QOL |
|---|---|---|---|---|---|---|---|
| 800 | 20 | LL | L | 900 | L | 700 | 200 |
| 500 | 40 | LD | L | 600 | D | * | * |
| 900 | 20 | DL | D | * | L | 800 | * |
| 300 | 20 | DD | D | * | D | * | * |

Table 9 (= Table 4 with Key Covariate observed)

| $X$ | Principal stratum | $Z_i$ | $S_i(T)$ | $\bar{Y}_i(T)$ | $S_i(C)$ | $\bar{Y}_i(C)$ |
|---|---|---|---|---|---|---|
| 800 | LL | T | L | 900 | ? | ? |
| 800 | LL | C | ? | ? | L | 700 |
| 500 | LD | T | L | 600 | ? | ? |
| 500 | LD | C | ? | ? | D | * |
| 900 | DL | T | D | * | ? | ? |
| 900 | DL | C | ? | ? | L | 800 |
| 300 | DD | T | D | * | ? | ? |
| 300 | DD | C | ? | ? | D | * |

distribution of the never jointly observed potential outcomes under $T$ and under $C$.

Again, our example will be extreme to illustrate ideas, and the actual methods of analysis with real examples will nearly always involve methods of analysis based on EM or MCMC methods for mixture models (e.g., Dempster, Laird and Rubin, 1977; Titterington, Smith and Makov, 1985; Aitkin and Rubin, 1985). Specifically, suppose the four marginal distributions of the QOL potential outcomes within each of three principal strata where they are well defined are normal (Gaussian): $N(900, 70^2)$ and $N(700, 50^2)$ for $LL$ when treated and not, $N(600, 40^2)$ for $LD$ when treated, and $N(800, 60^2)$ for $LD$ when not treated. Suppose also that the investigators are confident that the distributions are normal, which could occur, for example, if the QOL scores were based on the average of a large set of test items about activities that the individuals can and cannot perform. The hypothetical variances tend to be larger under the active treatment because the drug has a nonadditive effect, being more effective for some than for others.

Nothing essential changes in Tables 1–3, except with the addition of the standard deviation associated with each mean. But Table 4, giving the observed data, is changed in an important way when the distributions are given. First, the treated group that lives is still observed to have mean 700, but its distribution is markedly nonnormal, with one-third having mean 900 and standard deviation 70 and two-thirds having mean 600 and standard deviation 40. These components are easily observable as different because of the assumed normality and the small within-component standard deviations. In more subtle situations, we would have to use far more sophisticated analysis methods, which generalize standard mixture modeling techniques cited earlier. Although, in this extreme example, we can distinguish between the one-third and two-thirds mixture components in the treated group that lives, we do not yet know which is $LL$ and which is $LD$, however. But we do know that $LL+LD$ comprise 60% of the population, because 60% live in the random half exposed to the active treatment; thus, either $LL$ is one-third of the 60%, that is, 20% of the population, with a $N(900, 70^2)$ distribution when treated, or is two-thirds of the 60%, that is, 40% of the population, with a $N(600, 40^2)$ distribution when treated.

Table 8 (= Table 2 with Key Covariate)
*Principal strata among the patients, each split by treatment assignment*

| $X$ | % population | Principal stratum | $Z_i$ | Treatment $S_i(T)$ | Treatment $\bar{Y}_i(T)$ | Control $S_i(C)$ | Control $\bar{Y}_i(C)$ | Treatment effect on QOL |
|---|---|---|---|---|---|---|---|---|
| 800 | 10 | LL | T | L | 900 | L | 700 | 200 |
| 800 | 10 | LL | C | L | 900 | L | 700 | 200 |
| 500 | 20 | LD | T | L | 600 | D | * | * |
| 500 | 20 | LD | C | L | 600 | D | * | * |
| 900 | 10 | DL | T | D | * | L | 800 | * |
| 900 | 10 | DL | C | D | * | L | 800 | * |
| 300 | 10 | DD | T | D | * | D | * | * |
| 300 | 10 | DD | C | D | * | D | * | * |

Moving on to the observed QOL distribution for the 40% who survive in the control group, we will observe a mean of 750 arising from a half/half mixture of $N(800, 60^2)$ and $N(700, 50^2)$, where one component represents $LL$ and one component represents $DL$, but which is which? Before addressing this question, we note that here, decomposing the two components within the surviving subjects in the control group is not as obvious as in the treated group, because in the control the means of the two components are only about one standard deviation apart—but this is a standard problem using the aforementioned mixture modeling algorithms. Now each component in the surviving control group is observed to be one-half of the 40% who live under control; thus, from the control data, both $LL$ and $DL$ comprise 20% of the population. But from the treatment group, we know that $LL$ is either 20% or 40% of the population; so combining both pieces of information, $LL$ must be 20%. Also, $LL$'s QOL distribution when treated must therefore be $N(900, 70^2)$ and $LD$'s QOL distribution when treated must be $N(600, 40^2)$.

Furthermore, when not actively treated, from the surviving control group, $LL$'s QOL distribution must be either $N(700, 50^2)$ or $N(800, 60^2)$—the observed data cannot distinguish these two possibilities because both $LL$ and $DL$ are exactly the same size, 20% of the population. Of course, if we knew which of $LL$ or $LD$ had a higher mean QOL under control or which had the larger variance under control, we would know which was $N(800, 60^2)$. This points out a fragility in the estimation: as the principal strata get closer to each other in size or closer in means and variances, the estimation becomes more difficult.

Nevertheless, we have recovered much of Table 1. What are uncertain are the mean values of the QOL column under control and the SACE: under control, the mean QOL is either 800 for $LL$ and 700 for $DL$, or 700 for $LL$ and 800 for $DL$; and the SACE is either +100 or +200; in either case treatment is preferable to control.

## 7. DISCUSSION

It is not surprising that if we have distributional information and good covariates, the estimation of the principal strata can be sharpened, even without assumptions such as stochastic ordering of the groups. Of course, in general, estimation must involve Markov chain Monte Carlo techniques, likelihoods generally will not be regular, and models may need to be assessed using posterior predictive checks (Rubin, 1984; Gelman, Meng and Stern, 1996; Gelman, Carlin, Stern and Rubin, 2004).

It is also important to realize that the joint conditional distribution of treatment potential outcomes $(S_i(1), Y_i(1))$ and control potential outcomes $(S_i(0), Y_i(0))$ given covariates is inestimable in the sense that the likelihood is free of parameters governing this joint conditional distribution (i.e., the posterior distribution of these parameters equals their prior distribution). To address this situation, sensitivity analyses (e.g., Rosenbaum and Rubin, 1983) and the creation of large-sample bounds (e.g., Manski, 2003; Zhang and Rubin, 2003, in this specific problem) could well be quite helpful and informative. This joint distribution can be relevant to an individual's decision-making for treatment versus control. Knowing the detail provided by Table 1 beyond that in Table 4, however, can be helpful for this decision, even without knowledge of the inestimable joint distribution. For example, a person may decide that he is in better health than the typical patient and is therefore more likely to have outcomes like those either in the $LL$ stratum, who are the healthiest group when treated, or in the $LD$ stratum, who are the healthiest group when not treated. Consequently, such a person would be particularly interested in learning why those in $DL$ die when treated (e.g., negative side effects versus enjoying a too vigorous life), and then use this information to make a more informed choice than directly available from Table 4.

In conclusion, I think that causal inference modeling using potential outcomes and principal stratification, with its explicit and transparent assumptions, has helped clarify situations that statisticians must confront when there exists censoring of outcomes due to death of units, and has led to the creation of an approach to estimation that can be quite beneficial in a variety of difficult settings across a variety of disciplines.